\documentclass[12pt]{article}
\usepackage{theorem,amsfonts,amssymb}
\newtheorem{theorem}{Theorem}[section]
\newtheorem{proposition}{Proposition}[section]
\newtheorem{lemma}{Lemma}[section]
\newtheorem{corollary}{Corollary}[section]
\theorembodyfont{\upshape}
\newtheorem{definition}{Definition}
\newtheorem{remark}{Remark}[section]
\newtheorem{example}{Example}[section]
\newtheorem{proof}{Proof}

\newtheorem{acknowledgement}{Acknowledgement}

\newcommand{\bt}{\begin{theorem}}
\newcommand{\et}{\end{theorem}}
\newcommand{\bl}{\begin{lemma}}
\newcommand{\el}{\end{lemma}}
\newcommand{\bp}{\begin{proposition}}
\newcommand{\ep}{\end{proposition}}
\newcommand{\bex}{\begin{example}}
\newcommand{\eex}{\end{example}}
\newcommand{\bc}{\begin{corollary}}
\newcommand{\ec}{\end{corollary}}
\newcommand{\bo}{\begin{proof}}
\newcommand{\eo}{\end{proof}}
\newcommand{\bd}{\begin{definition}}
\newcommand{\ed}{\end{definition}}
\newcommand{\br}{\begin{remark}}
\newcommand{\er}{\end{remark}}
\newcommand{\be}{\begin{enumerate}}
\newcommand{\ee}{\end{enumerate}}

\begin{document}

\title{Strong relative property $(T)$ and spectral gap of random walks}
\author{C. R. E. Raja}
\date{}
\maketitle

\let\ol=\overline
\let\epsi=\epsilon
\let\vepsi=\varepsilon
\let\lam=\lambda
\let\Lam=\Lambda
\let\ap=\alpha
\let\vp=\varphi
\let\ra=\rightarrow
\let\Ra=\Rightarrow
\let \Llra=\Longleftrightarrow
\let\Lla=\Longleftarrow
\let\lra=\longrightarrow
\let\Lra=\Longrightarrow
\let\ba=\beta
\let\ov=\overline
\let\ga=\gamma
\let\Ba=\Delta
\let\Ga=\Gamma
\let\Da=\Delta
\let\Oa=\Omega
\let\Lam=\Lambda
\let\un=\upsilon

\newcommand{\Ad}{{\rm Ad}}
\newcommand{\cK}{{\cal K}}
\newcommand{\pc}{T_{si}}

\newcommand{\Aut}{{\rm Aut}}
\newcommand{\cR}{{\cal R}}
\newcommand{\Spr}{{\rm Spr}}
\newcommand{\cG}{{\cal G}}
\newcommand{\cH}{{\cal H}}
\newcommand{\G}{{\mathbb G}}
\newcommand{\Z}{{\mathbb Z}}
\newcommand{\Q}{{\mathbb Q}}
\newcommand{\cN}{{\cal N}}
\newcommand{\cM}{{\cal M}}
\newcommand{\cS}{{\cal S}}
\newcommand{\cL}{{\cal L}}
\newcommand{\N}{{\mathbb N}}
\newcommand{\R}{{\mathbb R}}
\newcommand{\C}{{\mathbb C}}
\newcommand{\T}{{\mathbb T}}
\newcommand{\mP}{{\mathbb P}}

\begin{abstract}
We consider strong relative property $(T)$ for pairs $(\Ga , G)$
where $\Ga$ acts on $G$.  If $N$ is a connected Lie group and $\Ga$
is a group of automorphisms of $N$, we choose a finite index
subgroup $\Ga ^0$ of $\Ga$ and obtain that $(\Ga , [\Ga ^0, N])$ has
strong relative property $(T)$ provided Zariski-closure of $\Ga$ has
no compact factor of positive dimension.  We apply this to obtain
the following: $G$ is a connected Lie group with solvable radical
$R$ and a semisimple Levi subgroup $S$.  If $S_{nc}$ denotes the
product of noncompact simple factors of $S$ and $S_T$ denotes the
product of simple factors in $S_{nc}$ that have property $(T)$, then
we show that $(\Ga , R)$ has strong relative property $(T)$ for a
Zariski-dense closed subgroup of $S_{nc}$ if and only if
$R=[S_{nc},R]$.  The case when $N$ is a vector group is discussed
separately and some interesting results are proved.  We also
considered actions on solenoids $K$ and proved that if $\Ga$ acts on
a solenoid $K$, then $(\Ga ,K)$ has strong relative property $(T)$
under certain conditions on $\Ga$.  For actions on solenoids we
provided some alternatives in terms of amenability and strong
relative property $(T)$.  We also provide some applications to the
spectral gap of $\pi (\mu )=\int \pi (g) d\mu (g)$ where $\pi$ is a
certain unitary representation and $\mu$ is a probability measure.
\end{abstract}



\begin{section}{Introduction}

Let $G$ be a topological (Hausdorff) group and $\pi $ be a strongly
continuous unitary representation of $G$. We first recall the
following weak containment notion: we say that $\pi$ weakly contains
the trivial representation or we write $I\prec \pi$ if for each
compact set $K$ and for each $\epsi
>0$, there is a vector $v$ such that $\sup _{g\in K} || \pi (g)v
-v|| <\epsi ||v||$.  It is easy to see that existence of nontrivial
invariant vector implies weak containment. We would be looking at
the existence of nontrivial invariant vectors for various forms of
weak containment.  One such well known condition is property $(T)$:
we say that $G$ has property $T$ if $I\prec \pi$ implies $\pi$ has
nontrivial invariant vectors.  We mainly consider relative property
$(T)$ for triples $(G, H, N)$ and strong relative property $(T)$ for
pairs $(H, N)$.

A triple $(G, H ,N)$ consisting of a topological group $G$ and its
subgroups $H$ and $N$ is said to have relative property $(T)$ if for
any unitary representation $\pi$ of $G$ such that restriction of
$\pi$ to $H$ weakly contains the trivial representation of $H$, we
have $\pi (N)$ has nontrivial invariant vectors.

A pair $(H, N)$ consisting of topological groups $H$ and $N$ with
$H$ acting on $N$ by automorphisms is said to have strong relative
property $(T)$ if $(G, H, N)$ has relative property $(T)$ for $G=
H\ltimes N$: if $\Ga$ is a topological group acting on a topological
group $N$ by automorphisms, then the semidirect product $\Ga \ltimes
N$ is the product space $\Ga\times N$ with multiplication given by
$(\ap, x)(\ba , y) = (\ap \ba , x \ap (y))$ and $N$ (resp. $\Ga$) is
identified with the closed subgroup $\{(e,x)\mid x\in N \}$ (resp.
$\{(\ap , e)\mid \ap \in \Ga \}$) of $\Ga \ltimes N$ under the map
$x\mapsto (e,x)$ (resp. $\ap \mapsto (\ap ,e)$).

Strong relative property $(T)$ was considered by \cite{S1} to obtain
a characterization of algebraic groups with property $(T)$ and
related results.  Relative property $(T)$ was crucial in determining
property $(T)$ for various type of groups (see \cite{BHV}).

If $G$ is a topological group and $H, N$ are subgroups of $G$, then
$[H, N]$ denotes the closed subgroup generated by $\{ aba^{-1}b^{-1}
\mid a\in H ~~{\rm and}~~ b\in N \}$.

We will be looking at strong relative property $(T)$ for actions on
connected Lie groups and for actions on solenoids: compact connected
abelian group of finite dimension will be called a solenoid.

We first recall that a group of automorphisms of a connected Lie
group $G$ may also be viewed as a group of linear transformation on
the Lie algebra $\cG$ of $G$ by identifying each automorphism with
its differential.  Let $\Ga$ be a group of automorphisms of $G$.
Then for any $\Ga$-invariant subgroup $N$ of $G$, let $[\Ga , N]$ be
the closed subgroup generated by $ \{ \ap (g)g^{-1} \mid g\in N, \ap
\in \Ga \}$.  Then $[\Ga , G]$ is a closed connected normal subgroup
of $G$ invariant under any automorphism $\ap$ normalizing $\Ga$.  It
can be proved that if $[\Ga ', G] \not = G$ for a subgroup $\Ga '$
of finite index in $\Ga$, then $(\Ga , G)$ will not have strong
relative property $(T)$ (cf. Lemma \ref{lsrp}).  Thus, we would have
to look at $[\Ga ', G]$ for any finite index subgroup $\Ga '$ of
$\Ga$.  But there could possibly be plenty of finite index subgroups
of $\Ga$. Considering Zariski-closure of $\Ga$ we now choose a
"smallest" finite index subgroup of $\Ga$ with which we obtain our
result.

Let $\Ga _Z$ be the Zariski-closure of $\Ga$ in $GL(\cG )$ and $\Ga
_Z^0$ be the Zariski-connected component of identity in $\Ga _Z$.
Then $\Ga _Z^0$ is a closed normal subgroup of finite index in $\Ga
_Z$. Let $\Ga ^0 = \Ga \cap \Ga _Z^0$. Then $\Ga ^0$ is a closed
normal subgroup of finite index in $\Ga$.  So, $[\Ga ^0, G]$ is
$\Ga$-invariant.  The action of $\Ga ^0$ satisfies an interesting
alternative that any $\Ga ^0$-orbit is either singleton or infinite
and it can be easily seen that this alternative holds good for
connected subgroups of $GL(\cG )$, but $\Ga ^0$ need not be
connected.

\br In general, $\Ga _Z$ and $\Ga _Z^0$ need not be subgroups of
automorphisms of the Lie group $G$, however if $G$ is a simply
connected group, then $\Ga _Z$ consists of automorphisms of $G$. As
we will be working with action of $\Ga$ on $\cG$, $\Ga _Z$ as well
as $\Ga _Z^0$ is more convenient.\er

We will now state our result for actions on connected nilpotent Lie
groups: proof of the following uses the well known criterion of
relative property $(T)$ in terms of invariant measures on projective
spaces and the results on invariant groups of measures on projective
spaces (cf. \cite{Da} and \cite{Ra}).

\bt\label{cnt} Let $N$ be a connected nilpotent Lie group and $\Ga$
be a group of automorphisms of $N$.  Suppose $\Ga _Z$ has no compact
factor of positive dimension.  Then $[\Ga _0, N]$ is the maximal
subgroup of $N$ such that $(\Ga , [\Ga ^0, N])$ has strong relative
property $(T)$. \et

Let $G$ be a connected Lie group with solvable radical $R$ and a
semisimple Levi subgroup $S$. Let $S_{nc}$ denote the product of
noncompact simple factors in $S$ and $S_T$ be the product of simple
factors in $S$ that have property $(T)$.  Cornulier \cite{Co}
introduced $T$-radical $R_T$ defined by $R_T = \ol{S_T[S_{nc}, R]}$
and proved that $(S_{nc}, R_T)$ has strong relative property $(T)$
(cf. Remark 3.3.7 of \cite{Co}).  We now obtain a Zariski-dense
subgroup version of this result in the following form.  Let $X$ be
the Lie algebra of $[S_{nc}, R]$ and $\rho \colon S_{nc} \to GL (X)$
be defined by for $g\in S_{nc}$, $\rho (g)$ is the differential of
the restriction to $R$ of the innerautomorphism given by $g$.

\bt\label{co} Suppose $\Ga$ is a closed subgroup of $S_{nc}$ such
that $\rho (\Ga)$ is Zariski-dense in $\rho (S_{nc})$. Then the
following are equivalent

\be
\item [(1)] $(\Ga , R)$ and $(\Ga S_T, \ol{S_TR})$ have
strong relative property $(T)$;

\item [(2)] $(G, \Ga , R)$ and $(G, \Ga S_T, \ol{S_TR})$ have
relative property $(T)$;

\item [(3)] $R= [S_{nc},R]$. \ee

In general, $(\Ga , [S_{nc}, R])$ and $(\Ga S_T, R_T)$ (resp. $(G,
\Ga , [S_{nc}, R])$ and $(G, \Ga S_T , R_T)$) have strong relative
property $(T)$ (resp. have relative property $(T)$). \et

One of the important step in the proof of \cite{Co} is strong
relative property $(T)$ for certain semisimple Lie group actions on
simply connected nilpotent Lie groups (cf. Proposition 1.5 of
\cite{Co}) but our Theorem \ref{cnt} is more general.  As we shall
see condition $(3)$ in \ref{co} implies that solvable group has to
be nilpotent (see Corollary \ref{slc}).  So, we pay more attention
to actions on connected nilpotent Lie groups.

To extend Theorem \ref{cnt} to actions on general Lie groups, it is
sufficient to consider the remaining case of actions on semisimple
Lie groups.  This is considered in the following (see Proposition
\ref{ss}) which shows that actions on connected nilpotent Lie groups
is the crucial case.

Our techniques used in the proof of Theorem \ref{cnt} can also be
used to prove strong relative property $(T)$ for actions on
solenoids: any finite-dimensional compact connected abelian group
will be called a solenoid.  Let $K$ be a solenoid and $\Ga$ be a
group of automorphisms of $K$. Denote by $\hat K$, the dual group of
characters on $K$. Then $\hat K$ is a torsion free abelian group of
finite rank and hence for some $n$, $\Z ^n \subset \hat K \subset \
\hat K\otimes \Q\simeq \Q ^n$: $n$ is the dimension of $K$.  Any
automorphism $\ap$ of $K$, gives an automorphism $\hat \ap$ of $\hat
K$.  Since $\Z ^n \subset \hat K$, $\hat \ap$ extends to a linear
map of $\Q ^n$.  Thus, any group of automorphisms of $K$ may be
realized as a subgroup of $GL_n(\Q )$.

For a prime number $p$, let $\Q_p$ be the $p$-adic field and
$\Q_\infty = \R$.  Let $\Ga _p$ be the $p$-adic-closure of $\Ga$ in
$GL_n (\Q _p )$ for finite $p$.  Let $\Ga _\infty$ be the group of
$\R$-points of the Zariski-closure of $\Ga$ in $GL_n (\C )$.  We
prove the following for actions on solenoids.

\bt\label{ast} Let $K$ be a solenoid and $\Ga$ be a group of
automorphisms of $K$ such that the dual action of $\Ga$ on $\hat K$
has no nontrivial invariant characters.  Suppose one of the
following condition is satisfied:

\be

\item[$(c_p)$] $\Ga _p$ is topologically generated by $p$-adic
one-parameter subgroups for some finite prime $p$;

\item[$(c_\infty)$] the real connected component of $\Ga _\infty$ has no
nontrivial compact factor and the dual action of $\Ga$ has no finite
orbit.
\ee

Then $(\Ga , K)$ has strong relative property $(T)$ \et

\bex We now give examples of subgroups of $GL_n(\Q )$ that satisfies
$(c_p)$ for some $p$.

For $p=\infty$, take $\Ga \subset GL_n(\Z )$ and $\Ga$ is
Zariski-dense in (a finite extension of) $SL_n(\C )$ (or more
generally in a noncompact simple Lie group).  Then $\Ga$ satisfies
$(c_\infty )$ but $\Ga$ does not satisfy $(c_p)$ for any prime $p$
as $\Ga _p$ is contained in the compact group $GL_n(\Z _p)$.

Fix a prime number $p$.  For simplicity we will restrict our
attention to $n=2$.  Consider the group $\Ga$ generated by
$$\{ \pmatrix { 1& {k\over p^i} \cr 0 &1} , \pmatrix { 1 & 0 \cr
{l\over p^j} &1 } \mid k, l, i,j \in \Z \}.$$  Then $\Ga _p$ is the
closed subgroup generated by the two $p$-adic one-parameter
unipotent groups $$\{ \pmatrix { 1& x \cr 0 &1} \mid x \in \Q _p
\},~~~ \{ \pmatrix { 1 & 0 \cr y &1 } \mid y \in \Q _p \}.$$  Thus,
$\Ga $ verifies $C_p$.  It is easy to see that $\Ga _q$ is a compact
subgroup of $GL_2 (\Z _q)$ for any prime $q\not = p$.  Thus, $\Ga$
does not verify $(c_q)$ for any prime $q\not = p$.

\eex

We now look at application of strong relative property $(T)$ to
ergodic theory of random walks.  Let $\mu$ be a regular Borel
probability measure on $G$ and $\pi$ be any unitary representation
of $G$.  Consider the $\mu$-average $\pi (\mu )$ of $\pi$ defined by
$$\pi (\mu )(v) = \int \pi (g)v d\mu (g)$$ for any vector $v$.  Since
$\mu$ is a probability measure, $||\pi _\mu || \leq 1$.  Let $\Spr
(\pi (\mu ))$ be the spectral radius of $\pi (\mu )$.  $\Spr (\pi
(\mu ) )<1$ has many interesting consequences to the ergodic
properties of the contraction $\pi (\mu )$ and its iterates $\pi
(\mu ^n)$ (see \cite{JRT} and \cite{LW}).  It is easy to see that
$\Spr (\pi (\mu) )<1$ implies $\pi$ does not weakly contains the
trivial representation. Results in \cite{S0} and Corollary 8 of
\cite{BG} prove that this necessary condition is also sufficient for
any adapted probability measure $\mu$ (that is, the closed subgroup
generated by the support of $\mu$ is the whole group) on connected
semisimple Lie groups with finite center that have no nontrivial
compact factors.  We use relative property $(T)$ to extend these
results to a large class of connected Lie groups.

\bc\label{src}  Let $S$ be a connected semisimple Lie group having
no nontrivial compact factors and $R$ be a locally compact group
such that $R$ acts on $S$ and $(S, R)$ has strong relative property
$(T)$.  Let $G$ be a factor group of $S\ltimes R$. Then

\be
\item for any unitary representation $\pi$ of $G$ with $I\not \prec
\pi$ and for any adapted probability measure $\mu$ on $G$, $\Spr
(\pi (\mu ) )<1$.

\item   In particular, if $G$ is a connected Lie group satisfying
$(3)$ of Theorem \ref{co} and $S=S_{nc}$, then $\Spr (\pi (\mu )
)<1$ for any adapted probability measure $\mu$ on $G$ and for any
unitary representation $\pi$ of $G$ with $I\not \prec \pi$. \ee\ec

Interesting and important kind of unitary representations are
obtained by considering $G$-spaces.  Suppose $X$ is a $G$-space with
$G$-invariant measure $m$.  Then there is a unitary representation
corresponding to the $G$-space $X$ given by
$$\pi (g) (f) (x)= f(g^{-1}x) $$ for $g \in G$, $x \in X$ and
$f \in L^2(X, m)$.  We are mainly interested in the unitary
representation $\pi _0$ of $G$ obtained by restricting $\pi (g)$ to
the subspace $L^2_0(X,m) = \{ f\in L^2(X, m) \mid \int f dm =0 \}$.

In this situation of representations arising from $G$-spaces, the
point-wise convergence of $\pi (\mu ^n) (f)$ is also an interesting
problem to consider not just for $f\in L^2$ but also for $f \in L^p$
for all $p\geq 1$.  Using methods in \cite{JRT}, we derive the
following on point-wise convergence for adapted and strictly
aperiodic $\mu$ (that is, the smallest closed normal subgroup a
coset of which contains the support of $\mu$ is $G$) and last part
uses the result of \cite{BC}.

\bc\label{pnc} Let $G$ be a connected Lie group satisfying $(3)$ of
Theorem \ref{co} with $S=S_{nc}$ and $\mu$ be an adapted and
strictly aperiodic probability measure on $G$. Let $X$ be a
$G$-space with $G$-invariant probability measure $m$ and $\pi _0$ be
the associated unitary representation of $G$ on $L^2_0(X,m)$.
Suppose $I\not \prec \pi _0$. Then
$$ \pi (\mu ^n) f(x) = \int f(g^{-1}x) d\mu ^n(g)$$ converges
$m$-a.e. for any $f\in L^p(X,m)$ and $1< p<\infty$ and the norm
convergence holds for $p=1$ also.

In particular, if $\Da$ is a lattice in $G$, then for $f \in
L^p(G/\Da)$ and $1< p<\infty$,
$$ \pi (\mu ^n)(f)(x)=\int f(g^{-1}x) d\mu ^n(g) $$ converges a.e.
\ec

\end{section}

\begin{section}{preliminaries}

Let $G$ be a topological (Hausdorff) group and $H$ be a subgroup of
$G$. For a unitary representation $\pi$ of $G$, $\cH _\pi$ denotes
the Hilbert space on which $\pi$ is defined and $\pi |_H$ denotes
the restriction of $\pi$ to $H$.

We first recall the following well-known weak-containment
properties.

\bd We say that a unitary representation $\pi$ of a topological
group $G$ weakly contains the trivial representation and we write
$I\prec \pi$ if for each compact set $K$ and for each $\epsi
>0$, there is a vector $v\in \cH _\pi$ such that $\sup _{g\in K} ||
\pi (g)v -v|| <\epsi ||v||$.  \ed

\br Suppose $G$ is a locally compact $\sigma$-compact group. Then
using the theory of positive definite functions, one can easily see
that $I\prec \pi$ if and only if there is a sequence $(v_n)$ of unit
vectors such that $||\pi (g)v_n-v_n||\to 0$ for all $g$ in $G$.\er

\bd We say that a locally compact group $G$ has property $(T)$ if
any unitary representation of $G$ that weakly contains the trivial
representation has nontrivial invariant vectors.\ed

Structure of Lie groups and algebraic groups having property $(T)$
is well understood resulting in rich class of groups having property
$(T)$ (cf. \cite{S1} and \cite{W1}).  We refer to \cite{BHV} for
details on groups having property $(T)$.  We now look at relativized
versions of property $(T)$.

\bd Let $G$ be a topological group with subgroups $H$ and $N$. We
say that $(G, H, N)$ has relative property $(T)$ if for any unitary
representation $\pi$ of $G$, $I\prec \pi | _H$ implies $\pi (N)$ has
nontrivial invariant vectors. \ed

\bd Let $H$ be a topological group acting on a topological group $N$
by automorphisms.  We say that the pair $(H, N)$ has strong relative
property $(T)$ if $(H\ltimes N, H, N)$ has relative property $(T)$.
\ed

Relative versions of property $(T)$ are used in constructing new
examples of groups with property $(T)$ and in obtaining spectral
radius of random walks (cf. \cite{BG}).  We first prove the
following useful elementary result on strong relative property
$(T)$.

\bl\label{lsrp} Let $H$ be a locally compact group acting on a
locally compact group $N$ by automorphisms and $M$ be a closed
subgroup of finite index in $H$.

\be

\item If $(H, N)$ has strong relative property $(T)$, then
$(M, N)$ also has strong relative property $(T)$.

\item If $[M,N]\not = N$, then $(H,N)$ does not have
strong relative property $(T)$. \ee\el

\bo Let $\pi $ be a unitary representation of $M\ltimes N$ such that
$I\prec \pi |_M$.  By considering a subgroup of $M$ we may assume
that $M$ is normal in $H$.  Consider the induced representation
$\sigma$ from $\pi$ to $H\ltimes N$. Since $H/M$ is finite, the
space on which $\sigma$ is defined consists of functions $f\colon
H\ltimes N \to \cH _\pi$ satisfying
$$f(hx) = \pi (h)f(x), ~~ h\in M\ltimes N, ~~x\in H\ltimes N$$
and $\sigma$ is defined by $$\sigma (g)f(x) = f(xg)$$ for all
$x,g\in H\ltimes N$.

Take a system of coset representatives $\{ Mx_i \}$ of $M$ in $H$.
Let $E$ be a compact subset of $H$ and $\epsi >0$.  Then there is a
compact subset $F$ of $M$ such that $E\subset \cup _iFx_i$. Since
$I\prec \pi _M$, there is a vector $v$ such that $\sup _{h\in K} ||
\pi (h)v - v|| <\epsi ||v||$ where $K= \cup x_iFx_i^{-1}F_1$ and
$F_1$ is a finite subset of $M$ such that $\{x_ix_j \}\subset \cup
F_1x_k$.  Define $f \colon H\ltimes N \to \cH _\pi$ by $f(hx_i) =
\pi (h) v$ for all $h \in M\ltimes N$ and all $i$. Then it can
easily be seen that $f\in \cH _\sigma$.

Take $g \in F$ and $ h \in M\ltimes N$, if $x_ix_j = h_{ij}x_k$ for
some $k$ with $h_{ij}\in F_1$, then $\sigma (gx_j) f(hx_i) = \pi
(hx_igx_i^{-1}h_{ij})v$. This shows that $$||\sigma (gx_j)f-f||\leq
\sup _{a\in K} || \pi (a)v - v|| $$ and hence $$ \sup _{g \in
E}||\sigma (g) f -f|| <\epsi ||v|| = \epsi ||f||$$ Thus, $I\prec
\sigma |_H$.  If $(H,N)$ has strong relative property $(T)$, $\sigma
(N)$ has nontrivial invariant vectors. Thus, there is a non-zero
function $f\colon H\ltimes N \to \cH _\pi$ such that $f(hx)= \pi
(h)f(x)$ and $f(xg)=f(x)$ for all $h\in M\ltimes N$, $x\in H\ltimes
N$ and $g\in N$.  Since $f$ is non-zero, $f(x_i)$ is a non-zero
vector in $\cH_\pi$ for some $i$. Thus, we get $f(x_i) = f(x_ig) =
\pi (x_igx_i^{-1})f(x_i)$ for all $g\in N$. This implies that
$f(x_i)$ is a non-zero vector invariant under $\pi (N)$.  This
proves that $(M, N)$ has strong relative property $(T)$.

If $[M,N]\not = N$, then in view of the first part, we may assume
that $[H, N]\not =N$.  Replacing $N$ by $N/[H,N]$ we may assume that
$H$ is trivial on $N$. This shows that $H$ is a normal subgroup
$H\ltimes N$, that is $H\ltimes N = H\times N$.  Let $\rho$ be a
nontrivial irreducible unitary representation of $N$. Define $\pi $
on $H\ltimes N=H\times N$ by $\pi (h,x) = \rho (x)$. Then $I\prec
\pi |_H=I$ and $\pi(N)$ has no nontrivial invariant vectors as $\pi
|_N =\rho$ is nontrivial and irreducible. \eo

We recall the following nondegeneracy conditions on measures.

\bd A probability measure $\mu$ on a locally compact group $G$ is
called adapted (resp. strictly aperiodic) if the support of $\mu$ is
not contained in a proper closed subgroup (resp. if the support of
$\mu$ is not contained in a coset of a proper closed normal
subgroup).\ed

\begin{subsection}{Groups without finite-dimensional representations}

Let $\mu$ be an adapted and strictly aperiodic probability measure
on a locally compact group $G$ and $\hat G$ be the equivalent
classes of irreducible unitary representations of $G$.  Considering
the unitary representation $\oplus _ {\pi \in \hat G \setminus \{ I
\}} \pi $, we see that $\sup _{\pi \in \hat G \setminus \{ I
\}}||\pi _\mu || <1$ if $||\rho (\mu )||<1 $ for unitary
representations $\rho$ without invariant vectors.  It follows from
Theorem 1 of \cite{BG} that for groups having property $(T)$ and
without nontrivial finite-dimensional unitary representations,
$||\pi (\mu )||<1 $ for any unitary representation $\pi$ without
invariant vectors. Thus for such groups, $\sup _{\pi \in \hat G
\setminus \{ I \}}||\pi _\mu || <1$ and hence point-wise convergence
holds (see also \cite{JRT}). Such groups are useful in building new
classes of groups on which the strong and point-wise convergences
hold (see \cite{BG}). We now attempt to study this type of groups.

We say that a Hausdorff topological group $G$ has $(NCF)$ (resp.
$(NFU)$) if $G$ has no nontrivial compact factor (resp. $G$ has no
nontrivial finite-dimensional unitary representation).  We first
observe the following:

\bl\label{br1} Let $G$ be a locally compact group having $(NFU)$.
Then we have $[G, G]=G$. \el

\bo Since nontrivial locally compact abelian groups have nontrivial
finite-dimensional unitary representations, we have $[G, G]=G$.\eo

As our main results are concerned about connected Lie groups, we
study Lie groups having $(NFU)$ or $(NCF)$.

As irreducible unitary representations of compact groups are
finite-dimensional, $G$ has $(NCF)$ if $G$ has $(NFU)$ but the
converse need not be true: the abelian group $\Q _p$ of $p$-adic
numbers has $(NCF)$ but all its irreducible unitary representations
are one-dimensional. We now prove the converse for connected Lie
groups.

\bl\label{ncfe} Let $G$ be a connected Lie group and $S$ be a
semisimple Levi subgroup of $G$. Then the following are equivalent:
\be
\item $G$ has $(NCF)$;

\item $G$ has $(NFU)$;

\item $[G, G]=G$, $S$ has $(NCF)$ and $G$ is the smallest
closed normal subgroup containing the semisimple Levi subgroup $S$.
\ee Moreover, if $G$ has $(NCF)$, then $G$ is either trivial or
nonamenable. \el

\bo Suppose $G$ has $(NCF)$.  Let $\rho$ be a finite-dimensional
unitary representation of $G$. Since $[G, G] = G$, we have $[\rho
(G), \rho (G)]= \rho (G)$.  Let $V$ be the real Lie subalgebra of
${\rm End}(\cH _\rho)$ such that $V$ is the Lie algebra of the Lie
subgroup $\rho (G)$ of $GL(\cH _\rho )$. Then $V$ coincides with its
commutator subalgebra and hence $V$ is an algebraic subalgebra (see
Corollary 3, section 6.2, Chapter 1 of \cite{OV} or Theorem 15,
Chapter III of \cite{Ch}).  Since $\rho (G)$ is a connected Lie
subgroup of $GL(\cH _\rho )$, $\rho (G)$ is a subgroup of finite
index in a real algebraic group.  This implies that $\rho (G)$ is a
closed subgroup of $GL(\cH _\rho )$. Since $\rho$ is unitary, $\rho
(G)$ is compact. Since $G$ has $(NCF)$, $\rho $ is the trivial
representation.

Assume that $G$ has $(NFU)$.  Since compact groups have nontrivial
finite-dimensional representations, and $S$ is a factor of $G$, $S$
has $(NCF)$.  Let $M_G$ be the smallest closed normal subgroup of
$G$ containing $S$.  Let $R$ be the solvable radical of $G$. Then
$G= SR$ and so $G/M_G \simeq R/R\cap M_G$. This implies that $G/M_G$
is a solvable group but by Lemma \ref{br1} we get that $[G, G]=G$,
hence $G=M_G$.

Assume that $[G, G]=G$, $S$ has $(NCF)$ and $G$ is the smallest
closed normal subgroup containing $S$.  Let $H$ be a closed normal
subgroup of $G$ such that $G/H$ is compact and let $\phi \colon G\to
G/H$ be the canonical projection.  Then $\phi (S)$ is a semisimple
Levi subgroup of $G/H$.  Since $G/H$ is compact, semisimple Levi
subgroup $\phi (S)$ of $G/H$ is also compact.  Since $S$ has no
compact factors, $\phi (S)$ is trivial. Thus, $H$ is a closed normal
subgroup containing $S$, hence $H=G$. This proves that $G$ has
$(NCF)$.

We now prove the second part.  Suppose $G$ has $(NCF)$ and amenable.
Then $S$ has $(NCF)$ and amenable, hence $S$ is trivial.  This shows
that $G$ is solvable.  Since $[G,G]=G$, $G$ is trivial.

\eo

\end{subsection}

\end{section}

\begin{section}{Actions on vector spaces}

We first prove the following lemma characterizing strong relative
property $(T)$ for linear actions, a part of the proof uses
projection valued measure method of Furstenberg as in \cite{Bu},
\cite{HV} and \cite{S1}.

\bl\label{lst} Let $V$ be a finite-dimensional vector space over a
local field of characteristic zero and $\Ga$ be a locally compact
$\sigma$-compact group of linear transformations on $V$.  Then the
following are equivalent:

\be
\item [(1)] $(\Ga , V)$ has strong relative property $(T)$;

\item [(2)] for any $\Ga$-invariant subspace $W$ of $V$ such that the
action of $\Ga$ on $V/W$ is contained in a compact extension of a
diagonalizable group, we have $W=V$. \ee \el

Before we proceed to the proof we fix the following notation:  for a
measure $\lam $ on a locally compact group $G$ and an automorphism
$\ap $ of $G$, define the measure $\ap (\lam)$ by $\ap (\lam )(B) =
\lam (\ap ^{-1}(B))$ for any Borel subset $B$ of $G$.

\bo It is easy to see that $(1)$ implies $(2)$.  We now prove the
converse.  Suppose $(\Ga ,V)$ does not have strong relative property
$(T)$. Then there is a unitary representation $\pi$ of $\Ga \ltimes
V$ on a Hilbert space $\cH$ such that $I\prec \pi |_\Ga$ and $\pi
(V)$ has no nontrivial invariant vectors.  Let $(v_n)$ be a sequence
of unit vectors such that $||\pi ( h) v_n -v_n||\to 0$ for all $h
\in \Ga$. Let $\hat V$ be the dual of $V$. Since dual of quotient
subspaces of $V$ correspond to subspaces of $\hat V$, it is
sufficient to find a nontrivial $\Ga$-invariant subspace $U$ of
$\hat V$ such that dual action of $\Ga$ on $U$ is contained in a
compact extension of a diagonalizable group.

Let $P$ be the projection valued measure associated to the direct
sum decomposition of $\pi |_V$. For any vector $v\in \cH$, let $\mu
_v(B) = ||P(B)v||^2$ for any Borel subset $B$ of $\hat V$. Then $\mu
_v$ is a non-negative measure on $\hat V$.  It is easy to verify
that $h\mu _v = \mu _{\pi(h)v}$ for $h\in \Ga$ and $v\in \cH$. Since
$\pi (V)$ has no nontrivial invariant vectors, $P(\{0\})$ is the
trivial projection, hence all $\mu _v$ have full measure in $\hat
V\setminus \{ 0 \}$.  Let $\mP (\hat V)$ be the projective space
associated to $\hat V$ and $\vp \colon \hat V\setminus \{0 \} \to
\mP (\hat V)$ be the canonical quotient map. Then any $\ap \in
GL(\hat V)$ defines a transformation $\ol \ap$ on $\mP (\hat V)$ by
$\ol \ap (\vp (v)) = \vp (\ap (v))$ for all $v\in \hat V\setminus \{
0 \}$. For simplicity we denote $\ol \ap$ also by $\ap$.  Then we
have $\ap \vp = \vp \ap$. Now, take $\lam _n =\vp (\mu _{v_n})$.
Then $(\lam _n)$ is a sequence of probability measures on $\mP (\hat
V)$. Since $\mP (\hat V)$ is compact, by passing to a subsequence,
we may assume that $\lam _n \to \lam$ in the weak* topology for some
probability measure $\lam $ on $\mP (\hat V)$. For $h \in \Ga $,
since $||\pi (h)v_n -v_n ||\to 0$, $h\lam _n -\lam _n \to 0$ in the
total variation norm. This implies since $\lam _n \to \lam$ that
$h\lam = \lam$ for all $h \in \Ga$.

Let $L$ be the smallest quasi-linear variety (that is, a finite
union of subspaces) of $\hat V$ such that $\vp (L\setminus \{ 0 \})$
contains the support of $\lam$. Define
$$N_L = \{ g\in GL(\hat V) \mid g(L) =L \}$$ and
$$ I_L = \{ g\in GL(\hat V) \mid g(x) =x ~{\rm for ~ all} ~x\in \pi
(L\setminus \{0 \}) \}.$$ Then $I_L$ and $N_L$ are algebraic groups
and $I_L$ is a normal subgroup of $N_L$ (see \cite{Da} and \cite{Fu}
for further details). Let $G_\lam = \{ g\in N_L \mid g\lam = \lam
\}$. Then since $\pi (L\setminus \{ 0\})$ contains the support of
$\lam$, $I_L\subset G_\lam$.  By Corollary 2.5 of \cite{Da} in the
archimedean case and Proposition 1 of \cite{Ra} in the
non-archimedean case, $G_\lam /I_L$ is compact (see also \cite{Fu}):
in the real case this implies that $G_\lam$ is a group of real
points of an algebraic group defined over reals.  Let $\Ga '$ be the
image of $\Ga$ in $GL(\hat V)$ under the dual action. Since $h\lam =
\lam$ for all $h \in \Ga$, $\Ga '\subset G_\lam$. Let $U$ be the
subspace spanned by $L$. Since $L$ is $\Ga$-invariant, $U$ is
$\Ga$-invariant. It follows from the definition of $I_L$ that $I_L$
restricted to $U$ is a group of diagonalizable transformations.
Since $G_\lam /I_L$ is compact and $\Ga '\subset G_\lam$, we get
that the dual action of $\Ga$ on $U$ is contained in a compact
extension of a diagonalizable group.  \eo

As a consequence of this criterion we have the following alternative
for closed irreducible subgroups of $GL(V)$.

\bc\label{a1} Let $\Ga$ be a closed subgroup of $GL(V)$.

\be

\item If $V$ is $\Ga$-irreducible, then $(\Ga, V)$ has strong
relative property $(T)$ or $\Ga$ has polynomial growth.

\item If dimension of $V$ is two, then $(\Ga ,V)$ has strong
relative property (T) or $\Ga$ is amenable.\ee\ec

\bo Assume $V$ is $\Ga$-irreducible.  If $(\Ga , V)$ does not have
strong relative property $(T)$, then by Lemma \ref{lst}, there is a
proper $\Ga$-invariant subspace $W$ of $V$ such that $\Ga$ action on
$V/W$ is contained in a compact extension of an abelian group. Since
$V$ is $\Ga$-irreducible, $W=\{ 0 \}$.  Thus, $\Ga$ is contained in
a compact extension of an abelian group.  Since $\Ga$ is closed in
$GL(V)$, $\Ga$ has polynomial growth.

Assume dimension of $V$ is two.  If $V$ is $\Ga$-irreducible, then
the result follows from the first part as groups of polynomial
growth are amenable.  If $V$ is not $\Ga$-irreducible, then since
dimension of $V$ is two, $\Ga$ is solvable, hence amenable.  \eo

We will now put some sufficient conditions on linear actions to
obtain the strong relative property $(T)$ of the corresponding pair.

\bp\label{rpa} Let $V$ be a finite-dimensional vector space over
$\R$ and $\Ga$ be a locally compact $\sigma$-compact group of linear
transformations on $V$ such that  the dual action of $\Ga$ on $\hat
V$ has no nonzero finite orbit.

\be

\item [(1)] Suppose $\Ga _Z$ has no compact factor of positive dimension.
Then $(\Ga , V)$ has strong relative property $(T)$.

\item [(2)] Suppose $\Ga $ has no finite dimensional unitary
representation. Then $(\Ga , V)$ has strong relative property $(T)$.
\ee \ep

\br It is easy to see that conditions on $\Ga $ as in $(1)$ and
$(2)$ are not a necessity for strong relative property $(T)$ but it
can easily be seen that having no nonzero finite orbit for the dual
action is a necessity. \er

\br If $\Ga \subset GL(V)$ is finitely generated, then $\Ga$ is
residually finite, hence $\Ga$ has nontrivial finite-dimensional
unitary representations.  Thus, (2) is not applicable for finitely
generated $\Ga$ but there are many finitely generated $\Ga$ having
$\Ga _Z$ with no compact factor of positive dimension (for instance
$\Ga = SL(n, \Z)$), hence $(1)$ is applicable (see the following
Example \ref{r1}). \er

\bo Let $W$ be a $\Ga$-invariant subspace of $V$ such that the
action of $\Ga$ on $V/W$ is contained in a compact extension of a
diagonalizable group.

We now prove (1).  Suppose $\Ga _Z$ has no compact factor of
positive dimension and the dual action of $\Ga$ has no nonzero
finite orbit. Let $\Ga _Z$ be the group of $\R$-points of the
Zariski-closure of $\Ga $. Then $W$ is invariant under $\Ga _Z$ and
$\Ga _Z$-action on $V/W$ is contained in a compact extension of a
diagonalizable group. Since $\Ga _Z$ has no compact factor of
positive dimension and nontrivial connected abelian Lie groups have
nontrivial compact factor, we get that $\Ga _Z$ is finite on $V/W$.
If $W\not = V$, then this implies that the dual action of $\Ga _Z$
has a nonzero finite orbit in $\hat V$. Since $\Ga \subset \Ga _Z$,
the dual action of $\Ga$ has a nonzero finite orbit. This is a
contradiction to the assumption that the dual action of $\Ga$ on
$\hat V$ has no nonzero finite orbit. Thus, $W=V$, hence (1) follows
from Lemma \ref{lst}.

We now prove (2).  Suppose $\Ga $ has no nontrivial
finite-dimensional unitary representation and the dual action of
$\Ga $ has no nonzero finite orbit.  Since nontrivial subgroup of
compact extension of abelian groups have nontrivial
finite-dimensional unitary representations, $\Ga$ is trivial on
$V/W$.  If $W\not = V$, then this implies that the dual action of
$\Ga $ has nontrivial invariant vectors in $\hat V$.  This is a
contradiction to the assumption that the dual action of $\Ga$ on
$\hat V$ has no nonzero finite orbit. Thus, $W=V$ and hence Lemma
\ref{lst} implies that $(\Ga ,V)$ has strong relative property
$(T)$. \eo

\bex\label{r1} (1) If $\Ga \subset GL_n(\R )$ is a Zariski-dense
subgroup in a connected non-compact simple Lie group, then $\Ga _Z$
has no compact factor of positive dimension: for instance, $\Ga$ is
a lattice in a non-compact simple Lie group such as $SU(n,1),
SO(p,q)$.

(2) If $\Ga \subset GL_n(\R )$ is such that no subgroup of finite
index in the group of $\R$-points of the Zarkiski-closure of $\Ga$
has invariant vectors, then it may be shown that $\Ga$ has no finite
orbit: for instance if $\Ga$ is a lattice in $SL_n(\R )$ or
$Sp_{2n}(\R )$, then $\Ga$ has no nonzero finite orbit.

(3) In particular, if $\Ga$ is a lattice in a connected non-compact
simple Lie subgroup of $GL_n(\R )$ that has no invariant vectors in
$\R ^n$, then $\Ga _Z$ has no compact factor of positive dimension
and $\Ga$ has no finite orbit. \eex

\end{section}

\begin{section}{Actions on connected Lie groups}

We first consider actions on connected nilpotent Lie groups. Let $N$
be a connected nilpotent Lie group.  Then denote $N_i= {[N,
N_{i-1}]}$ for $i\geq 1$ and $N_0 = N$.  It can easily be seen that
each $N_i$ is a closed connected characteristic subgroup of $N$ and
$N_{i-1}/N_i$ is abelian.

We first prove a Lemma on actions of Lie groups.

\bl\label{uni} Let $G$ be a connected Lie group and $\Ga$ be a group
of automorphisms of $G$ such that $\Ga _Z$ has no compact factors of
positive dimension.  Suppose $H$ and $N$ are a closed connected
$\Ga$-invariant subgroup of $G$ such that $\Ga ^0$ is trivial on
$G/H$ and on $H/N$.  Then $\Ga ^0$ is trivial on $G/N$. \el

\bo Let $\cG, \cH, \cN$ be Lie algebras of $G, H,N$ respectively.
Since $\Ga ^0$ acts trivially on $G/H$ and on $H/N$, $\Ga _0$ on
$\cG/\cN$ is contained in a unipotent subgroup.  Since $\Ga ^0$ is a
subgroup of finite index in $\Ga$, $\Ga$ on $\cG /\cN$ is contained
in a finite extension of a unipotent group.  Since $\Ga _Z$ has no
compact factor of positive dimension, $\Ga _Z^0 $ is trivial on
$\cG/\cN$.  This implies that $\Ga ^0$ is trivial on $\cG /\cN$.\eo

\bo $\!\!\!\!\!$ {\bf of Theorem \ref{cnt}}\ \  Let $H= [\Ga ^0,
N]$. Since $N$ is connected, $H$ is a closed connected normal
subgroup of $N$, hence $H$ is a Lie group.  Let $M=[\Ga ^0 , H]$.
Then $M$ is a closed connected normal subgroup of $H$.  Then $\Ga
^0$ is trivial on $N/H$ and on $H/M$.  By Lemma \ref{uni}, $\Ga ^0$
is trivial on $N/M$.  Since $H= [\Ga ^0, N]$, $H=M$.  Replacing $N$
by $H$ we may assume that $N= [\Ga ^0, N]$.

Let $\cN $ be the Lie algebra of $N$ and $\exp \colon \cN \to N$ be
the exponential map. Since $N$ is nilpotent, there is a $k\geq 0$
such that $N_k\not = \{e \}$ and $N_{k+1}= \{ e \}$.  We prove the
result by induction on the dimension of $N$.  Suppose $N$ is
abelian. Then $\exp$ is a homomorphism.  Let $\pi$ be a unitary
representation of $\Ga \ltimes N$ such that $I\prec \pi |_\Ga $.  We
now claim that $\pi (N)$ has nontrivial invariant vectors.  Let $\pi
_1 \colon \Ga \ltimes \cN$ be given by $\pi _1(\ap, v) = \pi (\ap ,
\exp (v))$. Then $\pi _1$ is a unitary representation such that
$I\prec \pi _1|_\Ga$. Let $V_0 = [\Ga ^0, \cN]$.  Then $V_0$ is a
$\Ga$-invariant subspace of $\cN$. If $\Ga _1$ is a normal subgroup
of finite index in $\Ga $, let $W =[\Ga _1 , V_0]$.  Then $W$ is a
$\Ga $-invariant vector space and $\Ga $ is finite on $V_0/W$.  This
implies that $\Ga _Z^0$ is trivial on $V_0/W$, hence $\Ga ^0$ is
trivial on $V_0/W$.  Since $\Ga ^0$ is trivial on $\cN/V_0$, by
Lemma \ref{uni}, $\Ga ^0$ is trivial on $\cN /W$.  Since $V_0 = [\Ga
^0, \cN ]$, $W=V_0$. This implies that the no finite index subgroup
of $\Ga $ acts trivially on a nontrivial quotient of $V_0$.  By
Proposition \ref{rpa}, we get that $(\Ga , V_0)$ has strong relative
property $(T)$.  Thus, $\pi _1(V_0) = \pi (\exp (V_0))$ has
nontrivial invariant vectors. It is easy to see that $\Ga ^0$ acts
trivially on $N/\ol {\exp (V_0)}$.  This implies that $N=[\Ga ^0,
N]= \ol {\exp (V_0)}$. Since $\pi (\exp (V_0))$ has nontrivial
invariant vectors, $\pi (N)$ has nontrivial invariant vectors.  This
proves that $(\Ga , N)$ has strong relative property $(T)$.

Suppose $N$ is not abelian.  Then $N_k$ is a nontrivial closed
connected $\Ga$-invariant central subgroup of $N$.  Let $\cN _k$ be
the Lie algebra of $N _k$ and $\exp _k$ is the exponential map of
$\cN _k$ onto $N_k$.  Let $V _2$ be a nontrivial $\Ga$-irreducible
subspace of $\cN _k$ and $A= \ol{\exp _k(V_2)}$.  Since $\exp _k$ is
a local diffeomorphism, $A$ is a nontrivial closed connected
$\Ga$-invariant central subgroup of $N$. Applying induction
hypothesis to $\Ga$-action on $N/A$, $(\Ga , N/A)$ has strong
relative property $(T)$.

If a finite index normal subgroup $\Ga _2$ of $\Ga$ acts trivially
on $V_2$, then $A$ is in the center of $\Ga _2 \ltimes N$. Since
$(\Ga , N/A)$ has strong relative property $(T)$, by Lemma
\ref{lsrp}, $(\Ga _2, N/A)$ has strong relative property $(T)$.
Since $A\subset [N, N]$, Proposition 3.1.3 of \cite{Co} implies that
$(\Ga _2, N)$ has strong relative property $(T)$ (see also Remark
3.1.7 of \cite{Co}). This implies that $(\Ga , N)$ has strong
relative property $(T)$.

If no finite index subgroup of $\Ga$ acts trivially on $V_2$.  Since
$V_2$ is $\Ga$-irreducible, the dual action of $\Ga$ on the dual of
$V_2$ has no finite orbit. Hence by (1) of Proposition \ref{rpa},
$(\Ga , V_2)$ has strong relative property $(T)$. This implies that
$(\Ga , A)$ has strong relative property $(T)$.  Since $(\Ga , N/A)$
has strong relative property $(T)$, we get that $(\Ga , N)$ has
strong relative property $(T)$. \eo

\bex We now look at few examples of $(\Ga , N)$ that are relevant to
Theorem \ref{cnt}.

(i) Take $N= \R^n$ and $\Ga = SL_n(R)$ where $R$ is a subring of
$\R$ (for instance, $R= \Z, \Q , \R$). Then it is easy to see that
$\Ga _Z= SL_n(\R )$ and $\Ga ^0=\Ga$.  Here $[\Ga , \R ^n]=\R ^n$.
By Theorem \ref{cnt}, $(\Ga , N)$ has strong relative property
$(T)$.

(ii) Take $N = \{ (t,s, r \mid t,s\in \R ^n, ~ r\in \R \}$ to be the
$(2n+1)$-dimensional Heisenberg group with multiplication given by
$$(t,s,r) (t, s' , r') = (t+t' , s+s' , r+r ' +<t,s'> )$$ for
$t,s,t' , s',\in \R ^n$, $r,r'\in \R$ and $\Ga$ be any Zariski-dense
subgroup of $SL_{2n}(\R )$.  Then $[\Ga , N]=N$ and hence by Theorem
\ref{cnt}, $(\Ga , N)$ has strong relative property $(T)$.

\eex

Using general theory of connected semisimple Lie groups, we have the
following which extends Proposition 1.5 of \cite{Co}.

\bc\label{slc} Let $N$ be a connected solvable Lie group and $S$ be
a connected semisimple Lie subgroup of automorphisms of $N$.  If $S$
has $(NCF)$ and $[S, N]=N$, then $N$ is a nilpotent group with its
maximal compact subgroup contained in $[N,N]$ and $(\Ga, N)$ has
strong relative property $(T)$ for any Zariski-dense closed subgroup
$\Ga$ of $S$.\ec

\bo We first claim that $N$ is nilpotent.  Let $G= S\ltimes N$. Then
$G$ is a connected Lie group and $N$ is its solvable radical.  Let
$\rho \colon G\to GL(\cG )$ be the adjoint representation of $G$
where $\cG$ is the Lie algebra of $G$.  Let $\tilde G$ be the
algebraic closure of $\rho (G)$ in $GL(\cG )$.  Then Chevalley
decomposition implies that $\tilde G = \tilde S TU$ where $\tilde S$
is a semisimple Levi subgroup, $T$ is an abelian group consisting of
semisimple elements and $U$ is the unipotent radical.  Also,
solvable radical of $\tilde G $ is $TU$ and $[\tilde S,T]=\{ e \}$.
This implies that $[\tilde S,TU]\subset U$.  Since $\rho (S)$ is a
semisimple subgroup of $\tilde G$, replacing $\tilde S$ by a
conjugate we may assume that $\rho (S)\subset \tilde S$.  Since
$\rho (N)$ is a connected solvable normal subgroup of $\rho (G)$,
$\rho (N)\subset TU$. Hence $[\rho (S), \rho (N)]\subset U$.  Since
$[S, N] =N$, $\rho (N)\subset U$, hence $\rho (N)$ is nilpotent.
Since kernel of $\rho$ is the center of $G$, $N$ is a nilpotent
group.

Let $\Ga$ be a Zariski-dense subgroup of $S$.  Then since $S$ is
connected semisimple Lie group, the connected component of $\Ga _Z$
is $S$, hence $\Ga ^0 =\Ga$.  Since $S$ has $(NCF)$, $\Ga _Z$ has no
compact factor of positive dimension.  Let $H=[\Ga , N]$.  Then $H$
is a closed connected $\Ga$-invariant normal subgroup of $N$. Since
$H$ is connected and $\Ga$ is Zariski-dense in $S$, $H$ is
$S$-invariant.  Since $H=[\Ga , N]$, $\Ga$ is trivial on $N/H$ and
hence $S$ is trivial on $N/H$.  Since $[S,N]=N$, $H=[\Ga , N]=N$. By
Theorem \ref{cnt}, $(\Ga ,N)$ has strong relative property $(T)$.

In order to prove maximal compact subgroup of $N$ is contained in
$[N,N]$.  We assume that $N$ is abelian and prove that $N$ has no
compact subgroup.  Let $L$ be a maximal compact subgroup of $N$.
Then $L$ is $S$-invariant. Since $L$ is a torus and $S$ is
connected, $S$ acts trivially on $L$.  Let $\cN$ and $\cL$ be the
Lie algebras of $N$ and $L$ respectively. Since $S$ is semisimple
and $\cL$ is $S$-invariant, there is a $S$-invariant vector space
$V$ such that $V\oplus \cL = \cN$.  Let $\exp$ be the exponential
map of $\cN$ onto $N$.  If $\{ \exp (tv) \mid t\in \R \}$ is
relatively compact for some $v\in V$, then $\exp (tv) \in L$ for all
$t\in \R$.  This implies that $v\in \cL$.  Since $V\cap \cL = \{ 0
\}$, $v=0$.  Thus, $\exp (V)$ is closed.  Since $V\oplus \cL =\cN$,
$S$ is trivial on $N/\exp (V)$, hence $[S, N]\subset \exp (V)$.
Since $[S ,N]=N$, we get that $N=\exp (V)$, hence $L$ is trivial.

\eo

\bex  We now give an example to show that $N$ in Corollary \ref{slc}
may have nontrivial compact subgroup.  Take $N = \{ (t,s, r+\Z )
\mid t,s, r \in \R \}$ to be the $3$-dimensional reduced Heisenberg
group with multiplication given by
$$(t,s,r+\Z) (t, s' , r'+\Z) = (t+t' , s+s' , r+r'+ts'+\Z )$$ for
$t,s,r ,t' , s',r'\in \R$ and $S=SL_{2}(\R )$.  Then $[S,N]=N$ and
$\{(0,0,r+\Z )\mid r\in \R \}$ is a compact central subgroup of $N$
of dimension one. \eex

We now prove Theorem \ref{co}.

\bo $\!\!\!\!\!$ {\bf of Theorem \ref{co}}\ \  Using the surjective
homomorphism $S\ltimes R \to G$ given by $(x,g)\mapsto gx$ we get
that $(1)$ implies $(2)$.  If $[S_{nc}, R]\not = R$, then let $N=
[S_{nc},R]$.  Then $N$ is a closed normal subgroup of $G$ and
$S_{nc}$ acts trivially on $R/N$.  This implies that $(\Ga, R)$ does
not have strong relative property $(T)$.  Since $G/S_{nc}R$ is
compact, $(G, \Ga , R)$ does not have relative property $(T)$. This
proves that $(2)$ implies $(3)$. Assume $[S_{nc}, R]=R$.  It follows
from Corollary \ref{slc} that $(\Ga , R)$ has strong relative
property $(T)$. Since $S_T\subset S_{nc}$ and $S_T$ has property
$(T)$, $(1)$ follows.

The second part may be proved from the first part by considering the
Lie group $S\ltimes R_T$ (instead of $S\ltimes R$) as $S_{nc}$ has
no nontrivial compact factor implies $[S_{nc}, R_T]=R_T$. \eo

\br We would like to remark that Theorem \ref{co} may be proved for
any locally compact $\sigma$-compact Zariski-dense (not necessarily
closed) subgroup $\Ga$ of $S_{nc}$.  \er

We now discuss extension of Theorem \ref{cnt} for actions on
connected Lie groups. It is quite clear that one needs to consider
the two distinct cases of actions on connected solvable Lie groups
and actions on connected semisimple Lie groups.  We now discuss the
case of actions on connected semisimple Lie groups.

\bp\label{ss} Let $\Ga$ be a group of automorphisms of a connected
semisimple Lie group $G$.  If $(\Ga, G)$ has strong relative
property $(T)$, then $G$ has property $(T)$. \ep

\bo  Suppose $(\Ga, G)$ has strong relative property $(T)$. Let $Z$
be the center of $G$. Then $Z$ is $\Ga$-invariant. It can easily be
seen that $(\Ga, G/Z)$ also has strong relative property $(T)$.
Using Theorem 1.7.11 of \cite{BHV} and replacing $G$ by $G/Z$, we
may assume that $G$ has no center and hence Aut$(G)$ is an almost
algebraic group (see \cite{Da1}). This implies that the connected
component of Aut$(G)$ has finite index in Aut$(G)$.  Since $G$ is a
connected semisimple Lie group, the group of inner automorphisms of
$G$ is the connected component of Aut$(G)$ (see Chapter III, Section
10.2, Corollary 2 of \cite{Bo}). Thus, there is a subgroup $\Ga _1$
of finite index in $\Ga$ such that $\Ga _1$ is a group of
inner-automorphisms on $G$, hence since $G$ has no center, $\Ga _1$
is isomorphic to a subgroup $G_1$ of $G$.  We may identify $\Ga _1$
with the subgroup $G_1$ of $G$.

Let $\pi$ be a unitary representation of $G$ such that $I\prec \pi$.
Define a unitary representation $\sigma $ of $\Ga _1\ltimes G$ by
$\sigma (x ,g) = \pi (gx)$ for all $(x,g)\in \Ga _1\ltimes G$. Then
$I\prec \sigma |_{\Ga _1}$.  Since $\Ga _1$ is a subgroup of finite
index in $\Ga$, by Lemma \ref{lsrp} we get that $(\Ga _1, G)$ has
strong relative property $(T)$ and hence $\sigma |_G= \pi$ has
nontrivial invariant vector. Thus, $G$ has property $(T)$. \eo

We now prove the spectral gap result.

\bo $\!\!\!\!\!$ {\bf of Corollary \ref{src}}\ \ Let $\vp \colon
S\ltimes R \to G$ be the canonical projection and $\pi$ be any
unitary representation of $G$. If $\mu$ is any adapted probability
measure on $G$ such that $\Spr (\pi (\mu ) )=1$. It follows from
\cite{BG} that $I\prec \pi \otimes \ol \pi$.  Let $\rho = \pi \circ
\vp$. Then $\rho$ is a unitary representation of $S\ltimes R$ and
$I\prec \rho \otimes \ol \rho$.  Since $S$ is a connected semisimple
Lie group having no nontrivial compact factors, Lemma 4 of
\cite{Ba1} implies that $I \prec \rho |_S$.  Since $(S, R)$ has
strong relative property $(T)$, $\rho (R)$ has nontrivial invariant
vectors. Since $R$ is normal, the space of $\rho (S)$-invariant
vectors is invariant and hence we may assume that $\rho (S)$ is
trivial. This proves that $I\prec \rho$ and hence $I\prec \pi$. \eo

\end{section}

\begin{section}{Solenoids}

We now look at actions on solenoids: recall that a compact connected
finite-dimensional abelian group is called solenoid.  Recall that if
$K$ is a solenoid (of dimension $n$), then $\Q _p^n$ may be realized
as a dense subgroup of $K$ and any group $\Ga$ of automorphisms can
be realized as a group of linear transformations on $\Q _p^n$.

\bp\label{cf} Let $K$ be a solenoid and $\Ga$ be a group of
automorphisms of $K$.  Suppose $(\Ga , K)$ does not have strong
relative property $(T)$.  Then we have the following:

\be

\item [(1)] for each $p$, there is a proper $\Ga$-invariant subspace $V_p$
of $\Q _p^n$ such that the action of $\Ga $ on $\Q _p^n /V_p$ is
contained in a compact extension of a diagonalizable group over $\Q_
p$;

\item [(2)] In addition if $\Ga$ is a finitely generated group and
for each $p$ either $\Ga _p$ is compact or the action of $\Ga $ is
irreducible on $\Q _p^n$, then there is a subgroup $\Ga _1$ of
finite index in $\Ga$ such that $\Ga _1$ has finite derived
subgroup. \ee \ep

\bo Let $p$ be a prime number or $p=\infty$.  Since $\hat K \subset
\Q ^n \subset \Q _p ^n$ and the dual of $\Q _p ^n$ is itself, we get
that there is a continuous homomorphisms $f_p \colon \Q _p ^n \to K$
such that $f_p(\Q _p ^n)$ is dense in $K$ and $f_p(\ap (x)) = \ap
(f_p(x))$ for all $\ap \in \Ga$ and $x\in K$.

Suppose $(\Ga , K)$ does not have strong relative property $(T)$.
Since $f_p(\Q _p ^n)$ is dense in $K$, $(\Ga , \Q _p ^n)$ also does
not have strong relative property $(T)$.  By Lemma \ref{lst}, there
is a proper $\Ga$-invariant subspace $W$ of $\Q _p ^n$ such that the
action of $\Ga$ on $\Q _p ^n /W$ is contained in a compact extension
of a diagonalizable group.  This proves (1).

We now prove (2).  Let $B$ be the dual of $\Q ^n$.  Then $K$ is a
quotient of $B$ and the $\Ga$ action on $K$ lifts to an action of
$B$.  Suppose $(\Ga , K)$ does not have strong relative property
$(T)$.  Then $(\Ga , B)$ also does not have strong relative property
$(T)$. Let $I$ be the set of all $p$ such that $\Ga _p$ is not
compact.  Since $\Ga$ action is irreducible on $\Q _p ^n$ or $\Ga
_p$ is compact, by Proposition \ref{lst} we get that $\Ga$ action on
$\Q _p ^n$ is contained in a compact extension of a diagonalizable
group, say $L_p$ for $p\in I$.  Since $\Ga$ is a finitely generated
group of matrices over $\Q$, $I$ is finite. Since $L_p$ is a compact
extension of a diagonalizable group, $L_p$ contains a subgroup
$L_p'$ of finite index such that $L_p'$ is a central group, in
particular derived group of $L_p'$ is compact. Let $\Ga _1 = \cap
_{p\in I} (L_p ' \cap \Ga)$. Since $I$ is finite, $\Ga _1$ is a
subgroup of finite index in $\Ga$. Since $L_p'$ has compact derived
group for any $p\in I$ and $\Ga _1$ is contained in $L_p'$, we get
that $\Ga _1$ has finite derived group. \eo

The next result shows the effectiveness of (2) of Proposition
\ref{cf} second result of which generalizes a result of Burger
\cite{Bu}: recall that action of a group $\Ga$ on a vector space is
called totally irreducible if the action of every subgroup of finite
index is irreducible.

\bc\label{irr}  Let $K$ be a solenoid and $\Ga$ be a finitely
generated group of automorphisms of $K$.

\be

\item If for each $p$ either $\Ga _p$ is compact or the action of $\Ga $
is irreducible on $\Q _p^n$,  then $(\Ga , K)$ has strong relative
property $(T)$ or $\Ga$ has polynomial growth.

\item If dimension of K is $2$, then $\Ga$ is amenable or $(\Ga, K)$
has strong relative property $(T)$.

\item If $\Ga \subset GL_n(\Z )$ and $\Ga$ is totally irreducible on
$\R^n$, then $(\Ga , K)$ has strong relative property $(T)$ or $K$
has dimension at most two.

\ee \ec

\br If $K$ is the $n$-dimensional torus, then $\Ga _p$ is compact
for any finite prime.  Thus, in this case, irreducibility conditions
in Corollary \ref{irr} needed only for the action of $\Ga$ on $\R
^n$. \er

\bo Assume that $(\Ga , K)$ does not have strong relative property
$(T)$. Suppose $\Ga$ is irreducible on $\Q _p ^n$ or $\Ga _p$ is
compact, (2) of Proposition \ref{cf} implies that $\Ga$ contains a
subgroup $\Ga _1$ of finite index such that $[\Ga _1, \Ga _1]$ is
finite. This in particular implies that $\Ga _1$ and hence $\Ga$
have polynomial growth. This proves the first part.

If dimension of $K$ is $2$ and $(\Ga , K)$ does not have relative
property $(T)$.  If $\Ga$ is not irreducible on some $\Q_p ^2$, then
$\Ga$ is contained in the group of upper triangular matrices which
is solvable, hence $\Ga$ is amenable.  So we may assume that $\Ga$
is irreducible on all $\Q_p ^2$.  By the first part, $\Ga$ is
amenable.  This proves the second part.

If $\Ga \subset GL_n(\Z )$ is totally irreducible on $\R^n$ and
$(\Ga , K)$ does not have relative property $(T)$.  Then $\Ga
\subset GL_n(\Z )$ implies that $\Ga _p$ is compact for any finite
prime $p$.  By the first part $\Ga$ is amenable.  By \cite{Ti},
$\Ga$ has a solvable subgroup $\Ga _1$ of finite index.  Since $\Ga$
is totally irreducible on $\R ^n$, $\Ga _1$ is irreducible on
$\R^n$. Since $\Ga _1$ is solvable, $n \leq 2$.

\eo

\bo $\!\!\!\!\!$ {\bf of Theorem \ref{ast}}\ \ Suppose $\Ga _p$ is
topologically generated by $p$-adic one-parameter subgroups.  If
$\Ga$ is contained in a compact extension of a diagonalizable group
over $\Q _p$, then let $H$ be a closed linear group over $\Q _p$ and
$D$ be a diagonalizable normal subgroup of $H$ such that $H/D$ is
compact and $\Ga$ is contained in $H$. Since there are no continuous
homomorphism from a $p$-adic one-parameter group into a compact
$p$-adic Lie group, we get that any $p$-adic one-parameter subgroup
of $H$ is contained in $D$. Since $D$ is diagonalizable, $H$ has no
$p$-adic one-parameter subgroup. Since $\Ga _p$ is generated by
one-parameter subgroups and $\Ga \subset H$, $\Ga $ is trivial.
Since dual action of $\Ga$ has no nontrivial invariant characters,
(1) of Proposition \ref{cf} is not satisfied.  Thus, $(\Ga, K)$ has
strong relative property $(T)$.

Suppose $\Ga$ satisfies $(c_\infty)$.  Let $\Ga _\infty ^0$ be the
connected component of identity in $\Ga _\infty$.  Then $\Ga _\infty
^0$ is a connected Lie group.  Let $S$ be the semisimple Levi
subgroup of $\Ga _\infty ^0$.  Since $S$ is a factor of $\Ga _\infty
^0$, $S$ is noncompact. This implies by Lemma \ref{ncfe} that $\Ga
_\infty ^0$ is either trivial or nonamenable. If $\Ga$ is a
contained in a compact extension of a diagonalizable group over
$\R$.  Since a compact extension of full diagonalizable group over
$\R$ is a $\R$-algebraic group, $\Ga _\infty$ is contained in a
compact extension of a diagonalizable group.  Since compact
extension of diagonalizable groups are amenable, we get that $\Ga
_\infty$ is amenable.  Since $\Ga _\infty ^0$ is either trivial or
nonamenable, we get that $\Ga _\infty ^0$ is trivial.  Since dual
action of $\Ga$ has no finite orbits, $(c_\infty)$ violates (1) of
Proposition \ref{cf}.  Thus, $(\Ga, K)$ has strong relative property
$(T)$.

\eo

We now prove spectral gap for actions on solenoids.

\bc\label{tfc} Let $K$ be a solenoid and $\Ga$ be a group of
automorphisms of $K$. Let $\pi _0$ be the unitary representation of
$\Ga$ on $L^2_0 (K) = \{ f\in L^2(K) \mid \int f = 0 \}$ given by
$\pi _0 (\ap )f (x) = f(\ap ^{-1}(x))$ for $\ap \in \Ga$ and $x\in
K$. Suppose $(\Ga , K)$ has strong relative property $(T)$.  Then
$\Spr (\pi _0 (\mu )) <1$ for any adapted probability measure $\mu$
on $\Ga$. \ec

\bo Let $\rho $ be the unitary representation of $\Ga \ltimes K$
given by $\rho (\ap ,a)(f) (x) = f(\ap ^{-1}(x)a)$ for all $(\ap ,
a) \in \Ga \ltimes K$, $x\in K$ and $f\in L^2_0 (K)$.  Then $\rho
|_\Ga = \pi _0$.  Also, any $f\in L^2_0(K)$ is $\rho (K)$-invariant
implies $f=0$.  Thus, $\rho (K)$ has no nontrivial invariant vectors
in $L^2_0 (K)$.

Suppose $(\Ga ,K)$ has strong relative property $(T)$.  Then $I\not
\prec \rho |_\Ga = \pi _0$ as $\rho (K)$ has no nontrivial invariant
vectors.  This implies since $\Ga$ is countable that $\Spr (\pi _0
(\mu )) <1$ for any adapted probability measure $\mu$ on $\Ga$ (cf.
\cite{LW}). \eo

\end{section}

\begin{section}{$G$-spaces}

\bo $\!\!\!\!\!$ {\bf of Corollary \ref{pnc}}\ \ Let $\mu$ be any
adapted and strictly aperiodic probability measure.  Then define
$\check \mu$ to be the probability measure defined by $\check \mu
(B) = \mu (B^{-1})$ for any Borel subset $B$ of $G$.  Let $\lam =
\sum _{n=1}^\infty {1\over 2^{n+1} }(\mu ^n *\check \mu ^n +\check
\mu ^n
*\mu ^n)$.  Then $\lam$ is a symmetric adapted probability measure on
$G$.

Let $\pi$ be a unitary representation of $G$ such that $I\not \prec
\pi$.  Then it follows from Corollary \ref{src} that $\Spr (\pi
(\lam ))<1$.  Since $\pi (\lam )$ is a self-adjoint positive
operator, $||\pi (\lam ) ||<1$.  This implies that $||\pi (\mu
)||<1$.

Let $X$ be any $G$-space with $G$-invariant measure $m$.  For any
$1\leq p < \infty$, let $\pi _p$ be the representation of $G$
defined on $L^p_0(X, m) = \{ f\in L^p(X, m) \mid \int f(x) dm(x) =0
\}$ given by $\pi _p(g)f(x) = f(g^{-1}x)$ for $g\in G$, $f\in
L^p_0(X, m)$ and $x\in X$.  Let $||\pi _p(\mu )||_p$ be the norm of
the operator $\pi _p(\mu )$ on $L^p_0(X, m)$.

Since $I\not \prec \pi _0=\pi _2$,  we get that $||\pi _2(\mu) ||_2
<1$. Since $||\pi _1(\mu )||_1 \leq 1$ and $||\pi _\infty (\mu
)(f)||_\infty \leq 1$ where $\pi _\infty$ is similarly defined on
$L^\infty$, we have by interpolation $||\pi _p (\mu )||_p < 1$ for
$1<p<\infty$ (cf. \cite{Ro}). This implies for $f\in L^p_0 (X, m)$
($1<p<\infty$), that $||\sum \pi _p(\mu ^n) f||_p \leq \sum ||\pi
_p(\mu ^n) f||_p <\infty$, hence $\lim \pi (\mu ^n) f(x) =0$ a.e.
Now the second part of the result follows from Theorem 1 of
\cite{BC} and from the first part. \eo
\end{section}

\begin{acknowledgement}
I would like to thank the referee for suggesting a version of
Theorem \ref{co} and for suggesting \cite{Co}. \end{acknowledgement}

\bigskip\medskip
\advance\baselineskip by 2pt
\begin{tabular}{ll}
C.\ R.\ E.\ Raja \\
Stat-Math Unit \\
Indian Statistical Institute (ISI) \\
8th Mile Mysore Road \\
Bangalore 560 059, India.\\
creraja@isibang.ac.in
\end{tabular}

\end{document}